\documentclass[10pt,twoside]{article}
\usepackage[latin1]{inputenc}
\usepackage{amsmath}
\usepackage{graphicx}
\usepackage{yhmath}
\usepackage[table]{xcolor}
\usepackage{mathrsfs} 
\usepackage{amssymb}
\usepackage{url}
\usepackage{makecell}
\usepackage{arydshln}
\usepackage{multirow}
\usepackage{arydshln}
\usepackage{mathtools}
\usepackage{stmaryrd}  

\input epsf
\def\limiten{\renewcommand{\arraystretch}{0.5}
\begin{array}[t]{c}\stackrel{}{\longrightarrow} \\
{\scriptstyle n\rightarrow
\infty}\end{array}\renewcommand{\arraystretch}{1}}

\def\limitepsn{\renewcommand{\arraystretch}{0.5}
\begin{array}[t]{c}\stackrel{a.s.}{\longrightarrow} \\
{\scriptstyle n \rightarrow
\infty}\end{array}\renewcommand{\arraystretch}{1}}

\def\limiteloin{\renewcommand{\arraystretch}{0.5}
\begin{array}[t]{c}\stackrel{{\cal D}}{\longrightarrow} \\
{\scriptstyle n\rightarrow
\infty}\end{array}\renewcommand{\arraystretch}{1}}

\def\limiteproban{\renewcommand{\arraystretch}{0.5}
\begin{array}[t]{c}\stackrel{{\cal P}}{\longrightarrow} \\
{\scriptstyle n\rightarrow
\infty}\end{array}\renewcommand{\arraystretch}{1}}

\def\equalpsn{\renewcommand{\arraystretch}{0.5}
\begin{array}[t]{c}\stackrel{a.s.}{=} \\
\end{array}\renewcommand{\arraystretch}{1}}

\numberwithin{equation}{section}

\newtheorem{thm}{Theorem}[section]

\newtheorem{lem}[thm]{Lemma}

\newtheorem{prop}[thm]{Proposition}

\newcommand{\E}{\ensuremath{\mathbb{E}}}
\newcommand{\R}{\ensuremath{\mathbb{R}}}
\newcommand{\Z}{\ensuremath{\mathbb{Z}}}

\newcommand{\N}{\ensuremath{\mathbb{N}}}

\newcommand{\prob}{\ensuremath{\mathbb{P}}}

\definecolor{grisclair}{gray}{0.9}

\renewcommand{\arraystretch}{.8}

\begin{document}
\title{\bf Epidemic change-point detection in general integer-valued time series}
 \maketitle \vspace{-1.0cm}
\begin{center}
   Mamadou Lamine DIOP \footnote{Supported by
   the MME-DII center of excellence (ANR-11-LABEX-0023-01) 
   } 
   and 
     William KENGNE \footnote{Developed within the ANR BREAKRISK: ANR-17-CE26-0001-01} 
     \footnote{This work was funded by CY Initiative of Excellence
(grant "Investissements d'Avenir" ANR-16-IDEX-0008),
Project "EcoDep" PSI-AAP2020-0000000013} 
 \end{center}

  \begin{center}
  { \it 
 THEMA, CY Cergy Paris Université, 33 Boulevard du Port, 95011 Cergy-Pontoise Cedex, France\\
  E-mail: mamadou-lamine.diop@u-cergy.fr ; william.kengne@u-cergy.fr  \\
  }
\end{center}

 \pagestyle{myheadings}
 \markboth{Epidemic change-point detection in general integer-valued time series}{Diop and Kengne}

~~\\
\textbf{Abstract}:
 In this paper, we consider the structural change in a class of discrete valued time series, which the true conditional distribution of the observations is assumed to be unknown. 
 The conditional mean of the process depends on a parameter $\theta^*$ which may change over time.
 We provide sufficient conditions for the consistency and the asymptotic normality of the Poisson quasi-maximum likelihood estimator (QMLE) of the model. 
 We consider an epidemic change-point detection and propose a test statistic based on the QMLE of the parameter.
%
Under the null hypothesis of a constant parameter (no change), the test statistic converges to a distribution obtained from a difference of two Brownian bridge.
The test statistic diverges to infinity under the epidemic alternative, which establishes that the proposed procedure is consistent in power.
The effectiveness of the proposed procedure is illustrated by simulated and real data examples.
 \medskip
 

\section{Introduction}
%
%
%
Change-point detection is a vast and active field of research, since its applications can be found in several areas, such as, epidemiology, finance, ecology, biology, etc. 
  This paper focuses on the epidemic change-point problem (see, for instance, \nocite{Levin1985} Levin and Kline (1985), \nocite{Yao1993} Yao (1993)) in a large class of integer-valued time series. The epidemic change-point problem involves testing the null hypothesis of no change versus the alternative that two changes occur during the data generating, with the structure of the first and the third segment is the same, and different from the second segment.

\medskip

  Assume that $\{Y_{t},\,t\in \Z \}$ is a time series of counts  
  and denote by $\mathcal{F}_{t-1}=\sigma\left\{Y_{t-1},\cdots ;X_{t-1},\cdots  \right\}$ the $\sigma$-field generated by the whole past at time $t-1$. 
Let $\Theta$ be a fixed compact subset of $\R^d$ ($d \in \N$) and $\mathcal T \subseteq \Z$. 
 For any $\theta \in \Theta$, define the class of integer-valued time series given by
 
 \medskip
 {\bf Class} $\mathcal{OD}_{\mathcal T }(f_{\theta})$: The process $Y=\{Y_{t},\,t\in \Z \}$ belongs to $\mathcal{OD}_{\mathcal T }(f_{\theta})$ if it satisfies:
  \begin{equation} \label{Model}
             \E(Y_t|\mathcal{F}_{t-1})=f_{\theta}(Y_{t-1},Y_{t-2},\cdots) ~~ \forall t \in \mathcal{T} ,
   \end{equation}
where $f_{\theta}(\cdot)$ is a measurable non-negative function, assumed to be known up to the  parameter $\theta$.  
Ahmad and Francq (2016) \nocite{Ahmad2016} carried out the inference question in the semiparametric setting in the class $\mathcal{OD}_{\Z }(f_{\theta})$, 
whereas Diop and Knegne (2020a, 2021) \nocite{Diop2020a} \nocite{Diop2021} focused on the model selection and multiple change-points problems in this class. 
Note that, numerous classical integer-valued time series models belong to the class  $\mathcal{OD}_{\Z }(f_{\theta})$: for instance, the Poisson INGARCH models (see for instance \nocite{Ferland2006} Ferland {\it et al.} (2006)), the negative binomial INGARCH models (proposed by \nocite{Zhu2011} Zhu (2011)), the binomial INGARCH (see \nocite{Weiss2014} Wei$ß$ and Pollett (2014)), the Poisson exponential autoregressive models (see \nocite{Fokianos2009} Fokianos {\it et al.} (2009)), 
the INAR models (see Wei$ß$ (2008) \nocite{Weiss2008} and  \nocite{Weiss2019} Wei$ß$ {\it et al.} (2019)).
%
%

\medskip

Firstly, we consider the class $\mathcal{OD}_{\Z}(f_{\theta^*})$ for some  $\theta^* \in \Theta$ and carry out the inference on the parameter $\theta^*$. The consistency and the asymptotic normality of the Poisson quasi-maximum likelihood estimator (QMLE) are addressed. These results are the same as the ones obtained by Ahmad and Francq (2016). But the conditions set here seem to be more straightforward than those needed by these authors. 
 
 \medskip

Secondly, we focus on the test with an epidemic alternative, for detecting changes in the parameter of the class $\mathcal{OD}_{\Z }(f_{\theta^*})$. The principle is that, the parameter has changed at time $t^*_1$, and then restored its original value after a time $t^*_2$; that is, an "epidemic" occurred between $t^*_1$ and $t^*_2$. 
This question has been addressed in several works; see, among others papers, Levin and Kline (1985), Yao (1993), \nocite{Csorgo1997} Cs{\"o}rg{\"o} and Horv{\'a}th (1997), \nocite{Ramanayake2003} Ramanayake and Gupta (2003), \nocite{Ravckauskas2004} Ra{\v{c}}kauskas and Suquet (2004), \nocite{Ravckauskas2006} Ra{\v{c}}kauskas and Suquet (2006), \nocite{Guan2007} Guan (2007), \nocite{Jaruvskova2011} Jaru{\v{s}}kov{\'a} and Piterbarg (2011), \nocite{Aston2012detecting, Aston2012Evaluating} Aston and Kirch (2012a, 2012b),  \nocite{Bucchia2014} Bucchia (2014), \nocite{Graiche2016} Graiche {\it et al.} (2016).
Most of these procedures are developed for epidemic change-point detection in the mean of random variables. Also, the case time series of count has not received great attention in the literature, while these models are very useful in many fields (see Section \ref{Sec_Real_data} for an example of application to the number of hospital admissions). 
For the general class $\mathcal{OD}_{\Z }(f_{\theta^*})$, we propose a test procedure based on the Poisson QMLE, for detecting epidemic change in the parameter $\theta^*$. 
Under the null hypothesis (no change) the test statistic converges to a distribution obtained from a difference between two Brownian bridge; this test statistic diverges to infinity under the epidemic alternative.  

\medskip

The paper is organized as follows. In Section 2, we set some assumptions, define the Poisson QMLE and establish its asymptotic properties. 
 Section 3 is devoted to the construction of the test statistic and the asymptotic studies under the null and the epidemic alternative. Some simulation results are displayed in Section 4.  Section 5 focuses on a real data example and Section 6 provides the proofs of the main results.

\section{Assumptions and Poisson QMLE}\label{Sect_Ass_PQMLE}
Throughout the sequel, the following notations will be used:
{\em
\begin{itemize}
 \item $ \|x \| \coloneqq  \sqrt{\sum_{i=1}^{p} |x_i|^2 } $, for any $x \in \mathbb{R}^{p}$;
%
%
 \item $ \|x \| \coloneqq  \sqrt{\sum_{i=1}^{p} \sum_{j=1}^{q} |x_{i,j}|^2 } $, for any matrix $x=(x_{i,j}) \in M_{p,q}(\R)$; where $M_{p,q}(\R)$ denotes the set of matrices of dimension $p\times q$ with coefficients in $\R$;
\item  $\left\|f\right\|_{\Theta} \coloneqq \sup_{\theta \in \Theta}\left(\left\|f(\theta)\right\|\right)$ for any function $f:\Theta \longrightarrow   M_{p,q}(\R)$;
\item $\left\|Y\right\|_r \coloneqq \E\left(\left\|Y\right\|^r\right)^{1/r}$, where $Y$ is a random vector with finite $r-$order moments; 
\item $T_{\ell,\ell'}=\{\ell,\ell+1,\cdots,\ell'\}$ for any $(\ell,\ell') \in \N^2$ such as $\ell \leq \ell'$.
\end{itemize}
}

 \medskip
    \noindent 
In the sequel, we will denote by $0$ the null vector of any vector space. 
Consider the following classical contraction condition on the function $f_\theta$.

    \medskip
    \noindent 
    \textbf{Assumption} \textbf{A}$_i (\Theta)$ ($i=0,1,2$):
    For any $y \in \mathbb{N}_0^{\N}$, the function $\theta \mapsto f_\theta(y)$ is $i$ times continuously differentiable on $\Theta$  with $ \left\| \partial^i f_\theta(0)/ \partial \theta^i\right\|_\Theta<\infty $; 
    and
      there exists a sequence of non-negative real numbers $(\alpha^{(i)}_k)_{k\geq 1} $ satisfying
     $ \sum\limits_{k=1}^{\infty} \alpha^{(0)}_k <1 $ (or $ \sum\limits_{k=1}^{\infty} \alpha^{(i)}_k <\infty $ for $i=1, 2$);
   such that for any  $y, y' \in \mathbb{N}_0^{\N}$,
  \[ 
  \sup_{\theta \in \Theta  } \Big \| \frac{\partial^i f_\theta(y)}{ \partial \theta^i}-\frac{\partial^i f_\theta(y')}{\partial\theta^i} \Big \|
  \leq  \sum\limits_{k=1}^{\infty}\alpha^{(i)}_k |y_k-y'_k|.
  \]

 \medskip 
  \noindent 
In the whole paper, it is assumed that any $\{Y_{t} ,\,t\in \Z\}$ belonging to $\mathcal{OD}_{\mathcal T }(f_{\theta})$ is a stationary and ergodic process satisfying: 
 \begin{equation}\label{moment}
    \exists C>0, \epsilon >0, \text{ such that } \forall t \in \Z, ~ ~ \E Y_{t}^{1+\epsilon} <C. 
   \end{equation}
   
   \medskip  
 \noindent 
Let $(Y_{1},\ldots,Y_{k}) \in \mathcal{OD}_{\{1,\ldots,k \}}(f_{\theta^{*}})$ be a trajectory with $k \geq 1$ and $\theta^* \in \Theta$. 
 Then, for any subset $\mathcal T \subseteq \{1,\cdots,k\}$, the conditional Poisson (quasi)log-likelihood computed on $\mathcal T$ is given (up to a constant) by
 \[
  L(\mathcal T, \theta) \coloneqq 
   \sum_{t \in \mathcal T} \ell_t(\theta) ~ \text{ with }~
  \ell_t(\theta) = Y_t\log \lambda_t(\theta)- \lambda_t(\theta),
  \]
  where $ \lambda_t(\theta) =f_\theta(Y_{t-1}, Y_{t-2}, \cdots )$.
 This conditional (quasi)log-likelihood is approximated (see also Ahmad and Francq (2016)) by
\begin{equation}\label{logvm}
\widehat{L}(\mathcal T, \theta) \coloneqq  \sum_{t\in \mathcal T} \widehat{\ell}_t(\theta) ~ \text{ with }~
 \widehat{\ell}_t(\theta) =  Y_t\log \widehat{\lambda}_t(\theta)- \widehat{\lambda}_t(\theta),
 \end{equation}
  where $ \widehat{\lambda}_t(\theta) =  f_\theta(Y_{t-1}, \cdots, Y_{1},0,\cdots)$. 
 According to (\ref{logvm}), the Poisson QMLE of $ \theta^*$ computed on $\mathcal T$ is defined by
 \begin{equation}\label{emv}
  \widehat{\theta}(\mathcal T) \coloneqq  \underset{\theta\in \Theta}{\text{argmax}} \big(\widehat{L}(\mathcal T,\theta)\big).
  \end{equation}

  \medskip
  
  \noindent
  When $(Y_{1},\ldots,Y_{n})$ is a trajectory of a process $\{Y_t,\, t \in \Z\}$ belonging to $\mathcal{OD}_{\Z}(f_{\theta^{*}})$, 
  we impose the following assumptions to study the asymptotic behavior of the Poisson QMLE.

  \begin{enumerate}
    
    \item [(\textbf{A0}):] 
     for all  $(\theta, \theta')\in \Theta^2$,
 $ \big( f_\theta(Y_{t-1}, Y_{t-2}, \cdots) \equalpsn  f_{\theta'}(Y_{t-1}, Y_{t-2}, \cdots)  ~ \text{ for some } t \in \N \big) \Rightarrow ~ \theta = \theta'$; 
 moreover, $\exists  \underline{c}>0$ such that $\displaystyle \inf_{ \theta \in \Theta} f_\theta(y)  \geq \underline{c}$, for all $ y \in  \N_0^{\N} $;

\item [(\textbf{A1}):] $\theta^* $ is an interior point of $\Theta \subset \mathbb{R}^{d}$;
 
  %
     \item [(\textbf{A2}):] for all $c \in \R^d$, $c' \frac{\partial \lambda_{t} (\theta^* )}{\partial    \theta} \equalpsn 0$   $\Rightarrow ~ c=0$.
    \end{enumerate}
    
    \medskip

\noindent 
The following proposition gives the strongly consistency and the asymptotic normality of the Poisson QMLE.
 \begin{prop}\label{Prop.1}
 Assume that $(Y_{1},\ldots,Y_{n})$ is a trajectory of a process $\{Y_t,\, t \in \Z\}$ belonging to $\mathcal{OD}_{\Z}(f_{\theta^{*}})$. 

\medskip
\rm
\noindent 
	(i.) \it If (\textbf{A0}), \textbf{A}$_0 (\Theta)$ and (\ref{moment}) (with $\epsilon \geq 1$) hold with 
	 \begin{equation}\label{Cond1.prpo1}
  \alpha_k^{(0)} = \mathcal O (k^{-\gamma}) ~~ for~ some~ \gamma>3/2,
 \end{equation}
 then 
 \[ \widehat{\theta}(T_{1,n}) \limitepsn  \theta^*.   \]
 \rm
 (ii.) \it If (\textbf{A0})-(\textbf{A2}), \textbf{A}$_i(\Theta)$ ($i=0,1,2$) and (\ref{moment}) (with $\epsilon \geq 3$) hold with 
	 \begin{equation}\label{Cond2.prpo1}
  \alpha_k^{(0)} +\alpha_k^{(1)}= \mathcal O (k^{-\gamma}) ~~ for~ some~ \gamma>3/2,
 \end{equation}
 then 
 \[ \sqrt{n}(\widehat{\theta}(T_{1,n})-\theta^*) \limiteloin \mathcal{N}(0,\Sigma)~ \text{ with } ~\Sigma \coloneqq J^{-1}_{\theta^*} I_{\theta^*} J^{-1}_{\theta^*}, \]
 where
 $J_{\theta^*} = \E \Big[ \frac{1}{\lambda_{0}(\theta^* )}  \frac{\partial \lambda_{0}(\theta^* )}{ \partial \theta} \frac{\partial \lambda_{0}(\theta^* )}{ \partial \theta'}  \Big] 
           ~~and~~ 
           I_{\theta^*} = \E \Big[\big( \frac{Y_0}{\lambda_{0}(\theta^* )}-1\big)^2  \frac{\partial \lambda_{0}(\theta^* )}{ \partial \theta} \frac{\partial \lambda_{0}(\theta^* )}{ \partial \theta'}  \Big]$.
 \end{prop}
This proposition will be proved by relying on some results which have already been established in Doukhan and Kengne (2015) without using  the  assumption of "conditional Poisson distribution".  
\medskip

\medskip
\noindent  
For any $\ell, \ell' \in \N$ with $\ell \leq \ell'$, define  the following matrices:
\begin{align*}
   &
   \widehat J(T_{\ell,\ell'})  =  \frac{1}{\ell'-\ell+1}\sum_{t \in T_{\ell,\ell'}}\frac{1}{\widehat \lambda_{t}(\widehat{\theta}(T_{\ell,\ell'}))}  \frac{\partial \widehat \lambda_{t}(\widehat{\theta}(T_{\ell,\ell'}))}{ \partial \theta} \frac{\partial \widehat \lambda_{t}(\widehat{\theta}(T_{\ell,\ell'}))}{ \partial \theta'} ,\\
   &
  \widehat I(T_{\ell,\ell'})  = \frac{1}{\ell'-\ell+1}\sum_{t \in T_{\ell,\ell'}} \Big(\frac{Y_t}{\widehat \lambda_{t}(\widehat{\theta}(T_{\ell,\ell'}))} -1 \Big)^2 \,  \frac{\partial \widehat \lambda_{t}(\widehat{\theta}(T_{\ell,\ell'}))}{ \partial \theta} \frac{\partial \widehat \lambda_{t}(\widehat{\theta}(T_{\ell,\ell'}))}{ \partial \theta'} .
 \end{align*}
According to (\textbf{A2}), one can easily show that the matrices $I$ and $J$ are symmetric and positive definite. Further, the part (i.) of Proposition \ref{Prop.1} implies the almost sure convergence of $\widehat J(T_{1,n})$ and  $\widehat I(T_{1,n})$ to $J_{\theta^*}$ and $I_{\theta^*}$, respectively. 
Therefore, $\widehat \Sigma_n  = \widehat J(T_{1,n}) \widehat I(T_{1,n})^{-1}   \widehat J(T_{1,n})$ is a consistent estimator of the covariance matrix $\Sigma $.


\section{Change-point test and asymptotic results}
Assume that $(Y_{1},\cdots,Y_{n})$ is an observed  trajectory  of the process $\{Y_{t},\,t\in \Z \}$  and 
we would like to test the null hypothesis of constant parameter 
\\
 \hspace*{.5cm}   
  $H_0$: $(Y_1,\cdots,Y_n)$  is a trajectory of the  process stationary $\{Y_{t},\,t\in \Z \} \in \mathcal{OD}_{\Z}(f_{\theta_1^*})$  with $\theta_1^* \in \Theta$,\\
 against the epidemic alternative   
\\  
\hspace*{.5cm}
$H_1$: there exists $(\theta^{*}_1,\theta^{*}_2,t^{*}_1,t^{*}_2) \in \Theta^{2}\times \{2,3,\cdots, n-1 \}^2$ (with $\theta^{*}_1 \neq \theta^{*}_2$ and $t^{*}_1 <t^{*}_2$) such that
     $(Y_1,\cdots,Y_{t^{*}_1})$ \\  \hspace*{1.2cm}
     and $(Y_{t^{*}_2+1},\cdots,Y_{n})$ are  trajectories of a process $\{Y^{(1)}_{t},\, t \in \Z\} \in \mathcal{OD}_{\Z}(f_{\theta^*_1})$, and $(Y_{t^{*}_1+1},\cdots,Y_{t^{*}_2})$ is a
     \\  \hspace*{1.2cm} trajectory of a process 
      $\{Y^{(2)}_{t}, \,t \in \Z\} \in \mathcal{OD}_{\Z}(f_{\theta^*_2})$.
\medskip

\medskip
\noindent  
 %
%
%
\noindent
We derive a retrospective test procedure in a semi-parametric setting, with a statistic based on the Poisson QMLE.
Suppose that $(u_n)_{n\geq 1}$ and $(v_n)_{n\geq 1}$ are two integer valued sequences such that:  $ u_n =o(n)$,  $v_n=o(n)$ and $ u_n,v_n \limiten +\infty$.
For all $n \geq 1$, define the matrix
\begin{multline*}\label{Sigma_un}
\widehat{\Sigma}(u_n)=\frac{1}{3} 
\big[
\widehat J(T_{1,u_n})  \widehat I(T_{1,u_n})^{-1}   \widehat J(T_{1,u_n}) 
+
\widehat J(T_{u_n+1,n-u_n})  \widehat I(T_{u_n+1,n-u_n})^{-1}  \widehat J(T_{u_n+1,n-u_n}) \\
+
\widehat J(T_{n-u_n+1,n})  \widehat I(T_{n-u_n+1,n})^{-1}  \widehat J(T_{n-u_n+1,n})
\big]
\end{multline*}
and the subset 
\[
\mathcal T_n = \left\{ (k_1, k_2) \in {([v_n, n-v_n] \cap \N)}^2 ~\text{ with }~  k_2- k_1  \geq v_n \right\}.
\]
For all $(k_1,k_2) \in \mathcal  T_n$, we introduce
 \begin{equation}\label{def_stat_C_n,k1,k2}
C_{n,k_1,k_2}= \frac{(k_2-k_1)}{n^{3/2}} \left[ \left(n-(k_2-k_1)\right)\widehat{\theta}(T_{k_1+1,k_2}) -  k_1 \widehat{\theta}(T_{1,k_1}) - (n-k_2)\widehat{\theta}(T_{k_2+1,n}) \right].
\end{equation}
Consider thus the test statistic given by
 \begin{equation}\label{Stat_Q_n}
 \widehat{Q}_n=\max_{(k_1, k_2) \in \mathcal T_n} \widehat{Q}_{n,k_1,k_2}
~\text{ with }~
 \widehat{Q}_{n,k_1,k_2}=C_{n,k_1,k_2}' \widehat{\Sigma}(u_n) C_{n,k_1,k_2}.
  \end{equation}

\noindent 
This statistic evaluates the distance between $\widehat{\theta}(T_{1,k_1})-\widehat{\theta}(T_{k_1+1,k_2}) $ and 
 $\widehat{\theta}(T_{k_2+1,n})-\widehat{\theta}(T_{k_1+1,k_2})$, for all $(k_1,k_2) \in \mathcal  T_n$.  These distances are not too large  
 in the absence of change-point (i.e., under $H_0$). Thus, the procedure rejects the null hypothesis if there exist two instants $k_1$ and $k_2$ such that the  distances exceed a suitably chosen constant.
 The subset $\mathcal T_n$ plays a very important role in the construction of the proposed test statistic, because it allows us to have  segments $T_ {1, k_1}$,  $T_ {k_1+1, k_2}$ and  $T_ {k_2+1, n}$ of sufficient large lengths, thus ensuring the convergence of the estimators computed.
The matrix $\widehat{\Sigma}(u_n)$  is also essentially useful to establish  the asymptotic properties of  $\widehat{Q}_n$, because:  
(i) under $H_0$,  each of the three matrices in the formula of $\widehat{\Sigma}(u_n)$ converges almost surely to the covariance matrix $\Sigma$
 and 
 (ii) under the epidemic alternative,  the first  and third matrices  converge to the covariance matrix of the stationary model of the first regime (or to the third regime)  which is positive definite.
The consistency of second matrix  is not ensured under the alternative; but it is positive semi-definite. 
Note that, a weight function can be used to increase the power of the test procedure based on the statistic $ \widehat{Q}_n$. See, for instance,   Doukhan and Kengne (2015) and Diop and Kengne (2017, 2020) for some examples. 
The statistic $\widehat{Q}_n$ can be seen as an extension to any parameter of the test statistic proposed by \nocite{Rackauskas2004} Rackauskas and Suquet (2004) (statistic $UI(n,\rho)$), \nocite{Jaruskova} Jarusková and Piterbarg (2011) (statistic $T_1^2$), \nocite{Bucchia2014} Bucchia (2014) (statistic $T_n(\alpha,\beta)$) or Aston and Kirch (2012a) (statistic $T_n^{B_2}$) in the context of mean change analysis. 
Indeed, in the particular case of the change-point detection in the mean with $\widehat{\theta}(\mathcal T)$ the empirical mean computed on the segment $\mathcal T$, the statistic $\widehat{Q}_n$ is equivalent to those proposed by these authors.  

\medskip
\noindent
 The following theorem establishes the asymptotic behavior of the test statistic under the null hypothesis. 
\begin{thm}\label{th1}
Under H$_0$ with $\theta^*_1 \in \overset{\circ}{\Theta}$, assume that (\textbf{A0})-(\textbf{A2}), \textbf{A}$_i(\Theta)$ ($i=0,1,2$) and (\ref{moment}) (with $\epsilon>3$)  and (\ref{Cond2.prpo1}) hold.
Then,
\begin{equation}\label{res_th1}
\widehat{Q}_n \limiteloin \sup_{0\leq \tau_1 < \tau_2 \leq 1}\left\|W_d(\tau_1)-W_d(\tau_2)\right\|^{2}, 
\end{equation}
where $W_d$ is a $d$-dimensional Brownian bridge.
\end{thm}

\noindent 
For a significance level $\alpha \in (0,1)$, the critical region of the test is then $(\widehat{C}_{n}>c_{d,\alpha})$, where $c_{d,\alpha}$ is the $(1-\alpha)$-quantile of the distribution
of $\sup_{0\leq\tau_1 < \tau_2\leq 1}\left\|W_d(\tau_1)-W_d(\tau_2)\right\|^{2}$.  This assures that the test procedure has correct size asymptotically. 
Table \ref{Table_C_alpha} below shows the values of $c_{d,\alpha}$ for $\alpha=0.01,\,0.05,\,0.10$ and $d=1,\ldots,5$, which are
obtained  by computing the empirical quantiles through Monte-Carlo simulations based on $5000$ replications.
The distribution was evaluated on a grid of size $1000$.

\begin{table}[h!]
\scriptsize
\centering
\caption{ \it 
Some empirical $(1-\alpha)$-quantiles of the distribution of $\underset{0\leq\tau_1 < \tau_2\leq 1}{\sup}\left\|W_d(\tau_1)-W_d(\tau_2)\right\|^{2}$.
}
\label{Table_C_alpha}
\vspace{.2cm}
\begin{tabular}{c c c c c c c c c c c}
\Xhline{.625pt}
 &&&&&\\
 &&\multicolumn{9} {l} {$d$}  \\
\cline{3-11}
 \rule[0cm]{0cm}{.25cm}
$\alpha$&&$1$&&$2$&&$3$&&$4$&&$5$\\
\Xhline{.6pt}
 \rule[0cm]{0cm}{.3cm} 
$0.01$&&$3.907$&&$7.320$&&$12.384$&&$16.004$&&$19.039$\\
 \rule[0cm]{0cm}{.3cm} 
$0.05$&&$2.973$&&$5.690$&&$8.948$&&$11.708$&&$14.471$\\
 \rule[0cm]{0cm}{.3cm} 
$0.10$&&$2.503$&&$4.988$&&$7.650$&&$9.954$&&$12.410$\\
 \Xhline{.625pt}
\end{tabular}
\end{table} 

\medskip
\noindent
Under the epidemic alternative, we set the following additional condition.

  \medskip
 {\bf Assumption B}: {\em There exists $(\tau^*_1,\tau^*_2) \in (0,1)^2$ such that $(t^*_1,t^*_2)=([n\tau^*_1], [n\tau^*_2])$ (with $[\cdot]$ is the integer part).}

\medskip

\noindent
Combining all the regularity assumptions given above, we obtain the following result.
\begin{thm}\label{th2}
Under $H_1$ with $\theta^*_1$ and $\theta^*_2$ belonging to  $\overset{\circ}{\Theta}$, assume that {\bf B}, (\textbf{A0})-(\textbf{A2}),  \textbf{A}$_i(\Theta)$ ($i=0,1,2$), (\ref{moment}) (with $\epsilon \geq 3$) and  (\ref{Cond2.prpo1}) hold. Then,

\begin{equation}\label{res_th2}
 \widehat{Q}_n  \limiteproban +\infty .
 \end{equation}
\end{thm}
 This theorem establishes the  consistency in power of the proposed procedure. 
 Under H$_1$, an estimator of the vector of breakpoints $\underline t^*=(t^*_1,t^*_2)$ is given by 
 \[   \widehat{ \underline t}_n =  \underset{(k_1, k_2) \in \mathcal T_n}{ \text{argmax}  }  C_{n,k_1,k_2}' \widehat{\Sigma}(u_n) C_{n,k_1,k_2} . \]


\section{Simulation study} 

We present some simulation results in order to assess the empirical size and power of the proposed  test procedure.
To do so, we consider the following processes:
\begin{flalign} 
&\bullet ~\text{\emph{Poisson-INGARCH processes}:} \nonumber\\
&\hspace{2.5cm} 
 \label{Poisson_model}  Y_{t}|\mathcal{F}_{t-1} \sim Poisson(\lambda_{t})~~\text{ with }~~ \lambda_{t}= \omega^*+ \alpha^* Y_{t-1} +\beta^* \lambda_{t-1},~ \text{for all }t \in \Z;\\
 &\bullet ~ \text{\emph{NB-INGARCH processes}:} \nonumber\\
 &\hspace{2.5cm}
   \label{NB_model} Y_{t}|\mathcal{F}_{t-1} \sim NB(r,p_{t})~~\text{ with }~~  r\frac{(1-p_{t})}{p_{t}}= \lambda_{t}= \omega^*+ \alpha^* Y_{t-1}+ \beta^* \lambda_{t-1}, ~ \text{for all }t \in \Z,&&
\end{flalign}
%
where  $NB(r, p)$ denotes the negative binomial distribution with parameters $r$ (assumed to be known) and $p$, and the parameter vector associated to the models is denoted by $\theta^* = (\omega^*,\alpha^*,\beta^*)$ which becomes $\theta^* = (\omega^*,\alpha^*)$ when $\beta^*=0$ (i.e, for an INARCH$(1)$ representation).
The NB-INGARCH processes are generated with $r = 5$.

\begin{figure}[h!]
\begin{center}
\includegraphics[height=10cm, width=17cm]{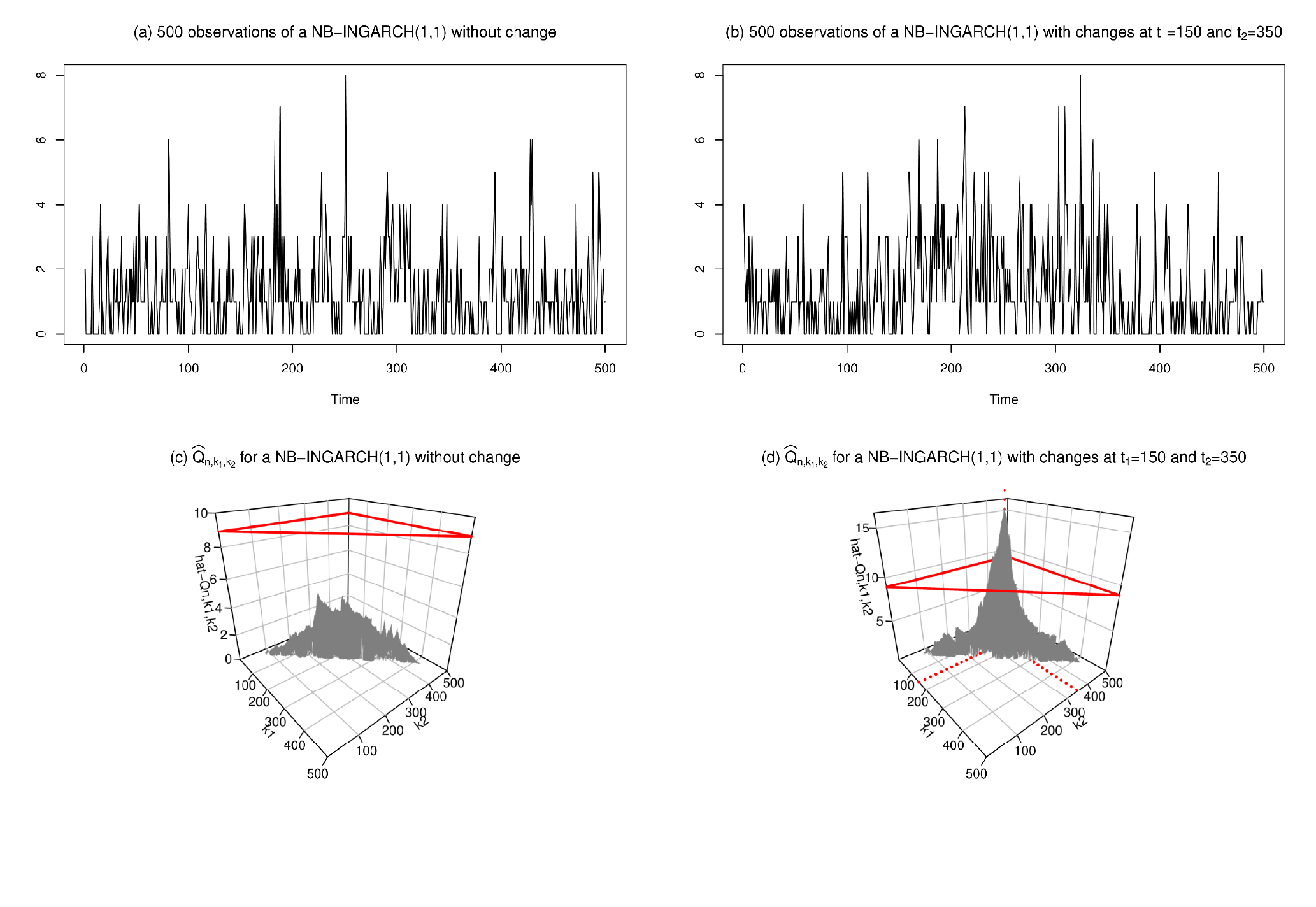} 
\end{center}
\vspace{-1.5cm}
\caption{\small \it Typical realization of 500 observations of two NB-INGARCH(1,1) processes with $r=5$ and the corresponding
statistics $\widehat{Q}_{n,k_1,k_2}$ for the epidemic change-point detection. 
(a) is a trajectory without change, where the true parameter $\theta^*=(0.5,0.2,0.35)$ is constant.
(b) is a trajectory generated under the epidemic alternative, where the parameter $\theta^*$ changes to $(1,0.2,0.35)$ at $t_1=150$ and  reverts back to $(0.5,0.2,0.35)$ at $t_2=350$. 
The horizontal triangles in (c) and (d) represent the limit of the critical region of the test, whereas 
the dotted lines show  the point where the maximum of $\widehat{Q}_{n,k_1,k_2}$ is reached.}
\label{Graphe_sim_NB}
\end{figure}

\medskip
Firstly, we generate two trajectories $(Y_1, \ldots, Y_{500})$ from (\ref{NB_model}): 
  a trajectory under $H_0$ with $\theta^*=(0.5,0.2,0.35)$ and a trajectory under $H_1$ with breaks at $t_1 = 150$ when $\theta^*$ changes to $(1,0.2,0.35)$ and $t_2 = 350$ when $\theta^*$ reverts back to $(0.5,0.2,0.35)$. 
  The procedure has been implemented on the R software (developed by the CRAN project).
  Figure \ref{Graphe_sim_NB} shows the realizations of the statistic $\widehat{Q}_{n,k_1,k_2}$ computed  with  $u_n= v_n=[\left(\log(n)\right)^{5/2}]$.  
As can be seen from this figure, in the scenario without change, the statistic $\widehat{Q}_{n,k_1,k_2}$ is less than the limit of the critical region that is represented by the horizontal triangle (see Figure \ref{Graphe_sim_NB}(c)). Under the alternative (of epidemic change),  $\widehat{Q}_{n,k_1,k_2}$ is greater than the critical value of the test and it reaches its maximum around the point where the changes occur (see the dotted lines in Figure \ref{Graphe_sim_NB}(d)).

\medskip
Now, for each of the two models (\ref{Poisson_model}) and (\ref{NB_model}), we are going to generate independent replications with sample size $n=500,\, 1000$ in the following situations: 
 a scenario where the parameter $\theta^*=\theta_0$ is constant (no change) 
 and a scenario where the parameter $\theta^*$ changes from $\theta_0$  to $\theta_1$ at time $t^*_1=0.3n$ and reverts back to $\theta_0$ at time $t^*_2=0.7n$. 
Table \ref{Table_res_sim} contains the empirical sizes and powers computed (under $H_0$ and $H_1$, respectively) as the proportion of the number of rejections of the null hypothesis based on $200$ repetitions. 
These results are obtained with a significance level $\alpha=5\%$.
 The scenario "$\theta_0=(22.75,0.18);~\theta_1=(14.5, 0.05)$"  considered here is related and close to the fitted representation obtained from the real data example (see below). 
As expected, the performance is better for the Poisson-INGARCH processes  than in the NB-INGARCH processes,  but the test procedure works well in both cases (see Table \ref{Table_res_sim}). 
 It produces reasonable empirical levels which are close to the nominal one when $n = 1000$. Also, the empirical powers increase with the sample size and are close to 1 when $n=1000$; which is consistent with the results of Theorem \ref{th2}.   

 \begin{table}[h!]
\centering
\scriptsize
\caption{\it Empirical sizes and powers at the nominal level $0.05$ for the epidemic change-point detection in the models (\ref{Poisson_model}) and (\ref{NB_model}).}
\label{Table_res_sim}
\vspace{.2cm}
\hspace*{-.3cm}
\begin{tabular}{llllcc}
\hline
\rule[0cm]{0cm}{.5cm}
&&&&$n=500$&$n=1000$\\
\hline
\rule[0cm]{0cm}{.4cm}
\multirow{6}*{\footnotesize Poisson-INGARCH processes}
 &Empirical levels: &&&&\\

                  &&$\theta_0=(22.75,0.18) $       &&0.040&0.045\\
                 
& && &&\\
                  && $\theta_0=(0.15,0.3,0.2)$   &&0.060&0.055\\
                
&&&&&\\
\rule[0cm]{0cm}{.1cm}
&Empirical powers:  && &&\\
                   &&$\theta_0=(22.75,0.18);$&$\theta_1=(14.5, 0.05);$ &0.995&1.000\\

& &&&&\\

                  && $\theta_0=(0.15,0.3,0.2);$&$\theta_1=(0.15,0.3,0.6)  ;$ &0.985&1.000\\   
                  
          & &&&&\\
   \hdashline[3pt/3pt]  
   & &&&&\\               
\multirow{6}*{\footnotesize NB-INGARCH processes}
 &Empirical levels: &&&&\\

                  &&$\theta_0=(22.75,0.18)$     &&0.030&0.040\\
                 
& && &&\\
                  && $\theta_0=(0.5,0.2,0.35)$   &&0.075&0.060\\
                
&&&&&\\
\rule[0cm]{0cm}{.1cm}
&Empirical powers:  && &&\\
                   &&$\theta_0=(22.75,0.18);$&$\theta_1=(14.5, 0.05);$ &0.980&1.000\\

& &&&&\\

                  && $\theta_0=(0.5,0.2,0.35);$&$\theta_1=(1,0.2,0.35) ;$ &0.965&0.990\\

 \Xhline{.75pt}
\end{tabular}
\end{table} 
         

 \section{Real data example}\label{Sec_Real_data}         
\begin{figure}[h!]
\begin{center}
\includegraphics[height=6.5cm, width=17cm]{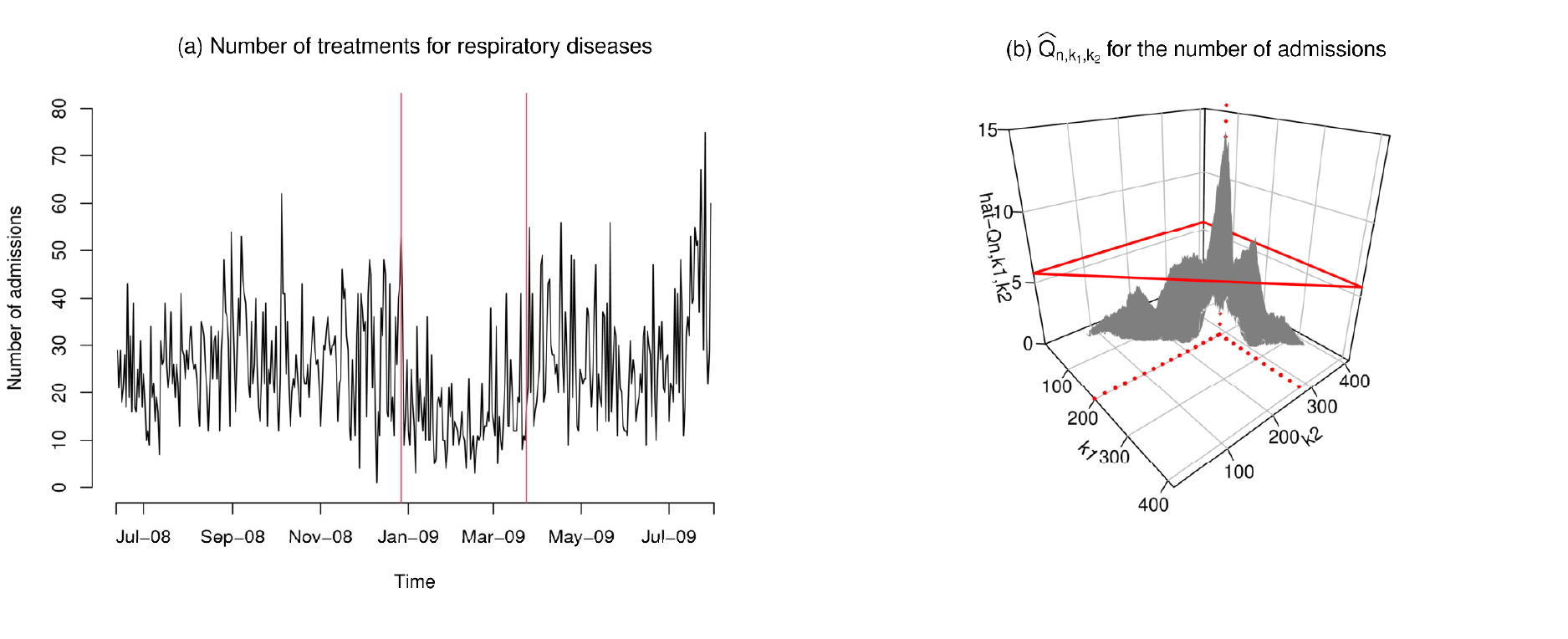} 
\end{center}
\vspace{-.8cm}
\caption{\small \it Plot of $\widehat{Q}_{n,k_1,k_2}$ for the epidemic change-point detection applied to 
the number of treatments for respiratory diseases in the Vit\'oria metropolitan area, Brazil, between June 13, 2008 and July 30, 2009 with an INARCH(1) representation.
The vertical lines in (a) are the estimated breakpoints. 
The horizontal triangle in (b) represents the limit of the critical region of the test, whereas 
the dotted lines show  the point where the maximum of $\widehat{Q}_{n,k_1,k_2}$ is reached.}
\label{Graphe_Application}
\end{figure}

 We investigate the  number of daily hospital admissions for respiratory diseases in children under $6$ years old in the Vit\'oria metropolitan area, Brazil.  The data are obtained from the Hospital Infantil Nossa Senhora da Gloria.  The time series is plotted in Figure \ref{Graphe_Application}(a).
 There are $413$ available observations that represent the admission from June 13, 2008 through July 30, 2009.   
 This time series is a part of a large dataset (available at https://rss.onlinelibrary.wiley.com/pb-assets/hub-assets/rss/Datasets/RSSC\%2067.2/C1239deSouza-1531120585220.zip) which has been studied by \nocite{Souza2018} Souza {\it et al.} (2018). 
 In their works, they used a hybrid generalized additive with Poisson marginal distribution to analyze the effects of some atmospheric pollutants on the number of hospital admissions due to cause-specific respiratory diseases.
 
The time series plot appears to show an epidemic change in the sequence.
To test this, we apply our detection procedure with an INARCH$(1)$ representation given by $\E(Y_t|\mathcal{F}_{t-1})= \lambda_{t}= \omega^*+ \alpha^* Y_{t-1}$. 
In each segment, to compute the QPMLE, the initial values $\lambda_{1}$ and $\partial \lambda_{1}/\partial \theta $ are set to be the empirical mean of the data and  the null vector, respectively.  
 For $u_n=[\left(\log(n)\right)^{5/2}]$ and $v_n=[\left(\log(n)\right)^{2}]$, Figure \ref{Graphe_Application}(b) shows the values of the statistic $\widehat{Q}_{n,k_1,k_2}$ corresponding to all the possible combinations $(k_1,k_2) \in \mathcal T_n$. 
 The critical value of the nominal level $\alpha=5\%$ is $c_{d,\alpha}=5.69$ and the resulting test statistic is $\widehat{Q}_{n}=14.72$; which implies the rejection of the null hypothesis (i.e., changes-points are detected).  
 The peak in the graph is reached at the point $(k_1,k_2)=(198, 285)$ which is the vector of the locations of the break-points estimated.  
The locations of the changes correspond to the dates December 27, 2008 and March 24, 2009. The second regime detected coincides with a large part of the austral summer which is from December to March; which partly explains the slight decrease of the number of hospital admissions observed in this period. 
The estimated model on each regime yields: 
 \begin{equation}\label{Estim_real_data}
\widehat \lambda_{t}=\left\{
\begin{array}{ll}
  \underset{(0.458)}{21.757} + \underset{(0.017)}{0.188}  Y_{t-1} ~\text{ for }~t\leq 198,  
  \\
\rule[0cm]{0cm}{.6cm}
\underset{(0.355)}{14.535} + \underset{(0.019)}{0.045}  Y_{t-1} ~\text{ for }~ 199 \leq t  \leq 285  , 
  \\
\rule[0cm]{0cm}{.6cm}
\underset{( 0.527)}{23.750} + \underset{(0.017)}{0.178}  Y_{t-1} ~\text{ for }~ t  \geq 286  , 
\end{array}
\right.
\end{equation}
where in parentheses are the robust standard errors of the estimators obtained from the sandwich matrix. 
In (\ref{Estim_real_data}), one can see that, the parameter of the first regime is very close to that of the third regime. This is in accordance with the alternative $H_1$ and lends a substantial support to the existence of an epidemic change-point in this series.

 \section{Proofs of the main results}  
  Let $ (\psi_n)_{n \in \N}  $ and $ (r_n)_{n \in \N}  $ be sequences of random variables or vectors. Throughout this section, we use the notation
  $ \psi_n = o_P(r_n)  $ to mean:  for all $  \varepsilon  > 0, ~ \prob( \|\psi_n \| \geq \varepsilon \|r_n \| ) \limiten 0$.
  Write $ \psi_n = O_P(r_n)$ to mean:  for all  $  \varepsilon > 0 $,  there exists $C>0$  such that  $\prob( \|\psi_n \| \geq C \|r_n \| )\leq \varepsilon $
   for $n$ large enough. 
   In the sequel, $C$ denotes a positive constant  whose the value may differ from one inequality to
another.

\subsection{Proof of Proposition \ref{Prop.1}}
 To simplify the expressions in this paragraph, we set: $\widehat{\theta}_n \coloneqq \widehat{\theta}_n(T_{1,n})$ and $L_n (\theta) \coloneqq L_n (T_{1,n},\theta)$ for any $\theta \in \Theta$.\\
  
	 (i.) To prove the first part of the proposition (i.e., the consistency), it suffices to show that the condition (4) of Ahmad and Francq (2016) is satisfied. This condition is established by  \nocite{Diop2020b} Diop and Kengne (2020b)  in their Remark 2.1 by using \textbf{A}$_0(\Theta)$ and (\ref{Cond1.prpo1}).\\
	
	(ii.) Applying the mean value theorem to the function $\theta \mapsto \frac{\partial}{\partial \theta_i} L_n (\theta)$  for all $i \in \{1,\cdots,d\}$, there exists 
	$\bar  \theta_{n,i}$ between $\widehat{\theta}_n$ and  $\theta^*$ such that
 \begin{equation*}\label{mvt_dl}
 \frac{\partial}{\partial \theta_i}L_n ( \widehat{\theta}_n)= \frac{\partial }{\partial \theta_i}L_n (\theta^*) + \frac{\partial^2}{\partial \theta \partial \theta_i}L_n (\bar \theta_{n,i}) (\widehat{\theta}_n-\theta^*),
\end{equation*}
which is equivalent to 
\begin{equation}\label{mvt.1_dl}
\sqrt{n}J_{\widehat \theta_{n}} (\widehat{\theta}_n-\theta^*) = 
\frac{1}{\sqrt{n}} \Big(\frac{\partial}{\partial \theta} L_n (\theta^*) -\frac{\partial}{\partial \theta}\widehat L_n ( \widehat{\theta}_n)\Big) 
+
\frac{1}{\sqrt{n}} \Big( \frac{\partial}{\partial \theta}\widehat L_n (\widehat{\theta}_n)-\frac{\partial}{\partial \theta}L_n (\widehat{\theta}_n)\Big)
\end{equation}
with
\begin{equation}\label{def_Jn}
J_{\widehat{\theta}_n} = \big(-\frac{1}{n}\frac{\partial^2}{\partial \theta \partial \theta_i} L_n (\bar \theta_{n,i})\big)_{1 \leq i \leq d}. 
\end{equation}
Moreover, by proceeding as in Lemma 7.1 of \nocite{Doukhan2015} Doukhan and Kengne (2015), we can show that
\[
\E\Big[ \frac{1}{\sqrt{n}} \Big\|\frac{\partial}{\partial \theta}\widehat L_n (\theta)- \frac{\partial}{\partial \theta} L_n (\theta)\Big\|_\Theta\Big]  \limiten 0.
\] 
In addition,  $\frac{\partial}{\partial \theta}  \widehat L_n ( \widehat{\theta}_n)=0$ for $n$ large enough since 
  $\widehat{\theta}_n$ is a local maximum of the function $\theta \mapsto  \widehat L_n (\theta)$ (from the assumption (\textbf{A1}) and the consistency of $\widehat{\theta}_n$).
  Thus, (\ref{mvt.1_dl}) gives 
 \begin{equation}\label{mvt.2_dl}
\sqrt{n}J_{\widehat \theta_{n}}  (\widehat{\theta}_n-\theta^*) = 
\frac{1}{\sqrt{n}} \frac{\partial}{\partial \theta} L_n (\theta^*) + o_P(1).
\end{equation} 
 The following lemma will be useful in the sequel.
 \begin{lem}\label{Lem00} 
 Assume that all the assumptions of Proposition \ref{Prop.1} hold. 
Then,
\begin{enumerate}
\rm
\item [(a)] 
$J_{\theta^*} = \E \Big[ \frac{1}{\lambda_{0}(\theta^* )}  \frac{\partial \lambda_{0}(\theta^* )}{ \partial \theta} \frac{\partial \lambda_{0}(\theta^* )}{ \partial \theta'}  \Big]<\infty 
           ~~and~~ 
           I_{\theta^*} = \E \Big[\big( \frac{Y_0}{\lambda_{0}(\theta^* )}-1\big)^2  \frac{\partial \lambda_{0}(\theta^* )}{ \partial \theta} \frac{\partial \lambda_{0}(\theta^* )}{ \partial \theta'}  \Big]<\infty .$
\rm
  \item [(b)]  \it $\left(\frac{\partial}{\partial \theta} \ell_t(\theta^*),\mathcal{F}_{t}\right)_{t \in \mathbb{Z}}$ is a stationary ergodic, square integrable martingale difference sequence with covariance matrix $I_{\theta^*}$;
    
\rm
\item [(c)]  \it $J_{\widehat \theta_{n}}  \limitepsn J_{\theta^*}$ and that the matrix $J_{\theta^*}$ is invertible. 
\end{enumerate}
 \end{lem} 
 {\bf Proof.}
 \begin{enumerate}
\item [(a)] It suffices to show that 
\[
(\text{a}_1):~ 
\E\Big[\Big\| \frac{1}{\lambda_{0}(\theta^* )}  \frac{\partial \lambda_{0}(\theta^* )}{ \partial \theta} \frac{\partial \lambda_{0}(\theta^* )}{ \partial \theta'}  \Big\|_\Theta\Big]<\infty 
~~\text{and}~~
(\text{a}_2):~
\E\Big[\Big\|\big( \frac{Y_0}{\lambda_{0}(\theta^* )}-1\big)^2  \frac{\partial \lambda_{0}(\theta^* )}{ \partial \theta} \frac{\partial \lambda_{0}(\theta^* )}{ \partial \theta'} \Big\|_\Theta\Big]<\infty.
\]
 $(\text{a}_1):$ From (\textbf{A0}), we have 
\begin{equation}\label{norm2_J}
\E\Big[\Big\| \frac{1}{\lambda_{0}(\theta^* )}  \frac{\partial \lambda_{0}(\theta^* )}{ \partial \theta} \frac{\partial \lambda_{0}(\theta^* )}{ \partial \theta'}  \Big\|_\Theta\Big]
 \leq 
\E\Big[\Big\| \frac{1}{\lambda_{0}(\theta^* )}\Big\|_\Theta \Big\| \frac{\partial \lambda_{0}(\theta^* )}{ \partial \theta}  \Big\|^2_\Theta\Big]
\leq 
\frac{1}{\underline{c}} \E\Big[\Big\| \frac{\partial \lambda_{0}(\theta^* )}{ \partial \theta}  \Big\|^2_\Theta\Big].
\end{equation}
Moreover, from Assumption \textbf{A}$_1 (\Theta)$, for all $t \in \Z$, we have
    \begin{align}\label{App_A1}
    \Big \|\frac{\partial \lambda_t(\theta)}{\partial \theta} \Big\|_\Theta
    &\leq  \Big\| \frac{\partial}{\partial \theta}f_\theta ( Y_{t-1},\cdots) - \frac{\partial}{\partial \theta}f_\theta  ( 0,\cdots) \Big\|_\Theta + \Big\|\frac{\partial}{\partial \theta}f_\theta  ( 0,\cdots) \Big\|_\Theta 
       \leq C +
       \sum\limits_{\ell \geq 1}  \alpha^{(1)}_{\ell}  Y_{t-\ell}.
      \end{align}
Then, applying the H\"{o}lder's inequality, we obtain 
\begin{align*}
\E\Big[\Big\| \frac{\partial \lambda_{t}(\theta^* )}{ \partial \theta}  \Big\|^2_\Theta\Big]
\leq
\Big( \Big\|\Big\|\frac{\partial \lambda_{t}(\theta^* )}{ \partial \theta}  \Big\|_\Theta\Big\|_4\Big)^2
 &\leq
\Big( \Big\|C+\sum\limits_{\ell \geq 1}  \alpha^{(1)}_{\ell}  Y_{t-\ell}\Big\|_4 \Big)^2\\
 &\leq \Big(C+\|Y_{0}\|_4\sum\limits_{\ell \geq 1}  \alpha^{(1)}_{\ell}  \Big)^2 <\infty 
~(\text{from } (\ref{moment}) \text{ with } \epsilon \geq 3). 
 \end{align*}
 We conclude the proof of $(\text{a}_1)$ from (\ref{norm2_J}). \\
 
 \noindent
 $(\text{a}_2):$ Again, using (\textbf{A0}), the H\"{o}lder's inequality and (\ref{App_A1}), we get 
 \begin{align*}
 \E\Big[\Big\|\big( \frac{Y_0}{\lambda_{0}(\theta^* )}-1\big)^2  \frac{\partial \lambda_{0}(\theta^* )}{ \partial \theta} \frac{\partial \lambda_{0}(\theta^* )}{ \partial \theta'} \Big\|_\Theta\Big]
 &\leq
 \frac{1}{\underline{c}}\left(\|Y_0\|_4+1\right)^2
 \Big( \Big\|\Big\|\frac{\partial \lambda_{0}(\theta^* )}{ \partial \theta}  \Big\|_\Theta\Big\|_4\Big)^2\\
 &\leq
 C\left(\|Y_0\|_4+1\right)^2
 \Big(C+\|Y_{0}\|_4\sum\limits_{\ell \geq 1}  \alpha^{(1)}_{\ell}  \Big)^2  <\infty ,
\end{align*} 
which achieves the proof of (a).
\item [(b)] See the proof of Theorem 2.2 of Ahmad and Francq (2016) or Lemma 7.2 of Diop and Kengne (2020b).
\item [(c)] This part is established in the proof of Theorem 3.2 of Doukhan and Kengne (2015) by using \textbf{A}$_i (\Theta)$ ($i=0,1,2$), (\ref{moment}) (with $\epsilon \geq 3$) and (\ref{Cond2.prpo1}). 
\end{enumerate}
\begin{flushright}
$\blacksquare$ 
\end{flushright}

 \noindent
Let us use the  Lemma \ref{Lem00} to complete the proof of the part (ii.) of Proposition \ref{Prop.1}. 
From Lemma \ref{Lem00}(c), for $n$ large enough such that  $J_{\bar{\theta}_n}$ (defined in (\ref{def_Jn})) is an invertible matrix. Then, the relation (\ref{mvt.2_dl}) is equivalent to
\begin{equation*}
\sqrt{n} (\widehat{\theta}_n-\theta^*) = 
J^{-1}_{\widehat \theta_{n}} \Big[ \frac{1}{\sqrt{n}} \frac{\partial}{\partial \theta_i} L_n (\theta^*) \Big] + o_P(1).
\end{equation*} 
Furthermore, applying the central limit theorem to the stationary ergodic martingale difference sequence $\left(\frac{\partial}{\partial \theta} \ell_t(\theta^*),\mathcal{F}_{t}\right)_{t \in \mathbb{Z}}$ (see  Lemma \ref{Lem00}(b)), we have 
\begin{equation*}
\frac{1}{\sqrt{n}} \frac{\partial }{\partial \theta} L_{n} (\theta^*) = \frac{1}{\sqrt{n}}\sum_{t=1}^{n}\frac{\partial}{\partial \theta} \ell_t(\theta^*) \limiteloin \mathcal{N}_d(0,I_{\theta^*}).
\end{equation*}
Therefore, for $n$ large enough, it holds that
\[
\sqrt{n} (\widehat{\theta}_n-\theta^*) = 
J^{-1}_{\theta^*} \Big[ \frac{1}{\sqrt{n}} \frac{\partial}{\partial \theta_i} L_n (\theta^*)\Big] + o_P(1)
\limiteloin \mathcal{N}_d (0,J^{-1}_{\theta^*}I_{\theta^*} J^{-1}_{\theta^*}).
\]
\begin{flushright}
$\blacksquare$
\end{flushright}
 
 %
%
 \subsection{Proof of Theorem \ref{th1}}
 \noindent 
 The following lemma is obtained from the Lemma A.1 and A.4 of Diop et Kengne (2021); the proof is then omitted. 
 %
 \begin{lem}\label{Lem0} 
 Assume that the assumptions of Theorem \ref{th1} hold. 
Then, 
  \[
(i.)~\frac{1}{n}\left\| \widehat{L}(T_{1,n},\theta)-L(T_{1,n},\theta)\right\|_{\Theta} \limitepsn 0
	~\text{ and }~
(ii.)~ \frac{1}{\sqrt{n}}\Big\| \frac{\partial \widehat{L}(T_{1,n},\theta)}{\partial \theta}  -\frac{\partial L(T_{1,n},\theta)}{\partial \theta}\Big\|_{\Theta}  \limiteproban 0.
\]
 \end{lem} 
  Define the statistic
  \[ 
  Q_n=\max_{(k_1, k_2) \in \mathcal T_n} Q_{n,k_1,k_2}
~\text{ with }~
 Q_{n,k_1,k_2}=C_{n,k_1,k_2}' \Sigma \,C_{n,k_1,k_2},
  \]
where $\Sigma$ is defined in Proposition \ref{Prop.1} and computed at $\theta^*_1$, under $H_0$.  
Consider the following lemma.
\begin{lem}\label{Lem1}
Assume that the assumptions of Theorem \ref{th1} hold. Then,
  \[
    \max_{(k_1,k_2) \in \mathcal T_n}\big|\widehat Q_{n,k_1,k_2}-Q_{n,k_1,k_2}\big|=o_P(1).
  \]
\end{lem}

\noindent
{\bf Proof.} \\
Let $(k_1,k_2) \in \mathcal T_n$. 
According to the asymptotic normality of the QMLE and the consistency of $\widehat{\Sigma}(u_n)$,  when $n \rightarrow \infty$, we have
\begin{equation}\label{cov_OP_H0}
   \left\{
\begin{array}{l} 
\big\| \widehat{\Sigma}(u_n)-\Sigma\big\|=o(1),~~~~ 
\big\|\sqrt{k_2-k_1} \left(\widehat{\theta}(T_{k_1+1,k_2}) -\theta^*_1\right)\big\|=O_P(1),\\
\rule[0cm]{0cm}{.6cm}

 \big\|\sqrt{k_1} \left(\widehat{\theta}(T_{1,k_1}) -\theta^*_1\right)\big\|=O_P(1), ~~~~
\big\|\sqrt{n-k_2} \left(\widehat{\theta}(T_{k_2+1,n}) -\theta^*_1\right)\big\|=O_P(1).
\end{array} 
\right.
\end{equation}
Then, we obtain
\begin{align*}
&\big|\widehat Q_{n,k_1,k_2}-Q_{n,k_1,k_2}\big|\\
&= \frac{(k_2-k_1)^2}{n^{3}}
 \Big|
\left[ \big(n-(k_2-k_1)\big)\widehat{\theta}(T_{k_1+1,k_2}) - \left( k_1 \widehat{\theta}(T_{1,k_1}) + (n-k_2)\widehat{\theta}(T_{k_2+1,n})\right) \right]'  \big( \widehat{\Sigma}(u_n)-\Sigma\big)\\
& \hspace{4.5cm} \times 
\left[ \big(n-(k_2-k_1)\big)\widehat{\theta}(T_{k_1+1,k_2}) - \left( k_1 \widehat{\theta}(T_{1,k_1}) + (n-k_2)\widehat{\theta}(T_{k_2+1,n})\right) \right]
\Big|\\
& \leq C\big\| \widehat{\Sigma}(u_n)-\Sigma\big\| \bigg[ 
\frac{(k_2-k_1)\big(n-(k_2-k_1)\big)^2 }{n^{3}} \big\|\sqrt{k_2-k_1} \left(\widehat{\theta}(T_{k_1+1,k_2}) -\theta^*_1\right)\big\|^2\\
& \hspace{1cm} +
\frac{k_1(k_2-k_1)^2 }{n^{3}} \big\|\sqrt{k_1} \left(\widehat{\theta}(T_{1,k_1}) -\theta^*_1\right)\big\|^2 + \frac{(n-k_2)(k_2-k_1)^2 }{n^{3}} \big\|\sqrt{n-k_2} \left(\widehat{\theta}(T_{k_2+1,n}) -\theta^*_1\right)\big\|^2
\bigg]\\
& \leq  o(1)O_P(1) = o_P(1).
\end{align*}
This allows to conclude the proof of the lemma.
\begin{flushright}
$\blacksquare$
\end{flushright}
Let $k,k^{\prime} \in[1,n]$, $\widetilde{\theta} \in \Theta$ and $i \in \{1,2,\cdots,d\}$. The mean value theorem to the function to $\theta \mapsto \frac{\partial}{\partial \theta_i} L({T_{k,k^{\prime}}},\theta)$ implies that there exists $\theta_{n,i}$ between $\widetilde{\theta}$ and $\theta^*_{1}$ such that
\[ 
\frac{\partial}{\partial \theta_i} L({T_{k,k^{\prime}}},\widetilde{\theta})=\frac{\partial}{\partial \theta_i} L({T_{k,k^{\prime}}},\theta^*_{1}) +\frac{\partial^{2}}{\partial \theta\partial \theta_i} L({T_{k,k^{\prime}}},\theta_{n,i})(\widetilde{\theta}-\theta^*_{1}).
\]
This is equivalent to
\begin{eqnarray}\label{Eq_Taylor}
(k^{\prime}-k+1)J_n({T_{k,k^{\prime}}},\widetilde{\theta}) (\widetilde{\theta}-\theta^*_{1})=\frac{\partial}{\partial \theta} L({T_{k,k^{\prime}}},\theta^*_{1})-\frac{\partial}{\partial \theta} L({T_{k,k^{\prime}}},\widetilde{\theta})
\end{eqnarray}
with
\begin{equation*}
J_n({T_{k,k^{\prime}}},\widetilde{\theta})= - \frac{1}{(k^{\prime}-k+1)}\frac{\partial^{2}}{\partial \theta\partial \theta_i} L({T_{k,k^{\prime}}},\theta_{n,i})_{1\leq i \leq d}.
\end{equation*}
%
%
%
%

\noindent  
We first use Lemma \ref{Lem00} to show that
\begin{equation}\label{Cond_proof_th1}
Q_n \limiteloin \sup_{0\leq \tau_1<\tau_2 \leq 1}\left\|W_d(\tau_1)-W_d(\tau_2)\right\|^{2}.
\end{equation}
Remark that 
\begin{align*}
C_{n,k_1,k_2}&= \frac{k_2-k_1}{n^{3/2}} \left[ (n-k_2)\left(\widehat{\theta}(T_{k_1+1,k_2})-\widehat{\theta}(T_{k_2+1,n})\right) -  k_1\left( \widehat{\theta}(T_{1,k_1}) -\widehat{\theta}(T_{k_1+1,k_2})\right) \right]\\
\text{and} \hspace{1.5cm}&\\
Q_{n,k_1,k_2}&=\big\|I^{-1/2}_{\theta^*_1}J_{\theta^*_1} C_{n,k_1,k_2}\big\|^2.
\end{align*}
 Let $(k_1,k_2) \in \mathcal T_n$. Applying (\ref{Eq_Taylor}) with $\widetilde{\theta}=\widehat{\theta}(T_{k_1+1,k_2})$ and $T_{k,k^\prime}=T_{k_1+1,k_2}$, we have
\begin{eqnarray}\label{Eq_T_(k_1+1,k_2)}
J_n(T_{k_1+1,k_2},\widehat{\theta}(T_{k_1+1,k_2})) \cdot (\widehat{\theta}(T_{k_1+1,k_2})-\theta_1^*)=\frac{1}{k_2-k_1}\Big(\frac{\partial}{\partial \theta} L(T_{k_1+1,k_2},\theta_1^*)-\frac{\partial}{\partial \theta} L(T_{k_1+1,k_2},\widehat{\theta}(T_{k_1+1,k_2}))\Big).
\end{eqnarray}
With $\widetilde{\theta}=\widehat{\theta}(T_{k_2+1,n})$ and $T_{k,k^\prime}=T_{k_2+1,n}$, (\ref{Eq_Taylor}) gives 
\begin{eqnarray}\label{Eq_T_(k_2+1,n)}
J_n(T_{k_2+1,n},\widehat{\theta}(T_{k_2+1,n})) \cdot (\widehat{\theta}(T_{k_2+1,n})-\theta_1^*)=\frac{1}{n-k_2}
\Big(\frac{\partial}{\partial \theta} L(T_{k_2+1,n},\theta_1^*)-\frac{\partial}{\partial \theta} L(T_{k_2+1,n},\widehat{\theta}(T_{k_2+1,n}))\Big).
\end{eqnarray}
Moreover, as $n\rightarrow +\infty$, Lemma \ref{Lem00}(c) (applied to $\theta^*_1$) implies
\begin{align*}
&\big\|J_n(T_{k_1+1,k_2},\widehat{\theta}(T_{k_1+1,k_2}))-J_{\theta^*_1}\big\|=o(1)
~\text{ and }~
\big\|J_n(T_{k_2+1,n},\widehat{\theta}(T_{k_2+1,n}))-J_{\theta^*_1}\big\|=o(1).
\end{align*}
Then, 
according to (\ref{cov_OP_H0}), for $n$ large enough, it holds from (\ref{Eq_T_(k_1+1,k_2)}) that 
\begin{align*}
\sqrt{k_2-k_1}J_{\theta^*_1}\left(\widehat{\theta}(T_{k_1+1,k_2})-\theta_1^*\right)
&=
\frac{1}{\sqrt{k_2-k_1}}\Big(\frac{\partial}{\partial \theta} L(T_{k_1+1,k_2},\theta_1^*)-\frac{\partial}{\partial \theta} L(T_{k_1+1,k_2},\widehat{\theta}(T_{k_1+1,k_2}))\Big) 
\\
& \hspace{1cm} -\sqrt{k_2-k_1}\big(\left(J_n(T_{k_1+1,k_2},\widehat{\theta}(T_{k_1+1,k_2}))-J_{\theta^*_1}\right)\left(\widehat{\theta}(T_{k_1+1,k_2})-\theta_0\right)\big)\\
&=
\frac{1}{\sqrt{k_2-k_1}}\Big(\frac{\partial}{\partial \theta} L(T_{k_1+1,k_2},\theta_1^*)-\frac{\partial}{\partial \theta} L(T_{k_1+1,k_2},\widehat{\theta}(T_{k_1+1,k_2}))\Big) +o_P(1)\\
&=
 \frac{1}{\sqrt{k_2-k_1}}\Big(\frac{\partial}{\partial \theta} L(T_{k_1+1,k_2},\theta_1^*)- \frac{\partial}{\partial \theta} \widehat L(T_{k_1+1,k_2},\widehat{\theta}(T_{k_1+1,k_2}))\Big)+o_P(1)\\
& \hspace{1cm} 
+\frac{1}{\sqrt{k_2-k_1}}\Big(\frac{\partial}{\partial \theta} \widehat L(T_{k_1+1,k_2},\widehat{\theta}(T_{k_1+1,k_2}))-\frac{\partial}{\partial \theta} L(T_{k_1+1,k_2},\widehat{\theta}(T_{k_1+1,k_2}))\Big) \\
&=
\frac{1}{\sqrt{k_2-k_1}}\Big(\frac{\partial}{\partial \theta} L(T_{k_1+1,k_2},\theta_1^*)-\frac{\partial}{\partial \theta} \widehat L(T_{k_1+1,k_2},\widehat{\theta}(T_{k_1+1,k_2}))\Big) +o_P(1),
\end{align*}
where the last equality is obtained from Lemma \ref{Lem0} ($ii.$). 
It is equivalent to
\begin{equation}\label{Eq_a}
J_{\theta^*_1}\left(\widehat{\theta}(T_{k_1+1,k_2})-\theta_1^* \right)=\frac{1}{k_2-k_1}\Big(\frac{\partial}{\partial \theta} L(T_{k_1+1,k_2},\theta_1^*)-\frac{\partial}{\partial \theta} \widehat L(T_{k_1+1,k_2},\widehat{\theta}(T_{k_1+1,k_2}))\Big) +o_P\Big(\frac{1}{\sqrt{k_2-k_1}}\Big).
\end{equation}
 For $n$ large enough, $\widehat{\theta}(T_{k_1+1,k_2})$ is an interior point of $\Theta$ and we have $ \frac{\partial}{\partial \theta} \widehat L(T_{k_1+1,k_2},\widehat{\theta}(T_{k_1+1,k_2}))=0$.
Thus, from (\ref{Eq_a}), we obtain
\begin{equation}\label{Eq_a_bis}
J_{\theta^*_1}\left(\widehat{\theta}(T_{k_1+1,k_2})-\theta_1^*\right)
 =\frac{1}{k_2-k_1}\frac{\partial}{\partial \theta} L(T_{k_1+1,k_2},\theta_1^*)+o_P\Big(\frac{1}{\sqrt{k_2-k_1}}\Big).
\end{equation}
Similarly, using (\ref{Eq_T_(k_2+1,n)}), we also obtain
\begin{equation}\label{Eq_a_bisbis}
J_{\theta^*_1}\left(\widehat{\theta}(T_{k_2+1,n})-\theta_1^*\right)=\frac{1}{n-k_2}\frac{\partial}{\partial \theta} L(T_{k_2+1,n},\theta_1^*)+o_P\Big(\frac{1}{\sqrt{n-k_2}}\Big).
\end{equation}
The subtraction of  (\ref{Eq_a_bis}) and (\ref{Eq_a_bisbis})  gives
\begin{multline*}
J_{\theta^*_1}\left(\widehat{\theta}(T_{k_1+1,k_2})-\widehat{\theta}(T_{k_2+1,n})\right)
=
\frac{1}{k_2-k_1}\frac{\partial}{\partial \theta} L(T_{k_1+1,k_2},\theta_1^*)-\frac{1}{n-k_2}\frac{\partial}{\partial \theta} L(T_{k_2+1,n},\theta_1^*)\\
+o_P\Big(\frac{1}{\sqrt{k_2-k_1}}+\frac{1}{\sqrt{n-k_2}}\Big);
\end{multline*}
i.e.,
\begin{multline}\label{part1_C_n,k1,k2}
\frac{(k_2-k_1)(n-k_2)}{n^{3/2}} J_{\theta^*_1}\left(\widehat{\theta}(T_{k_1+1,k_2})-\widehat{\theta}(T_{k_2+1,n})\right)=\\
\frac{1}{n^{3/2}}
\Big[
(n-k_2)\frac{\partial}{\partial \theta} L(T_{k_1+1,k_2},\theta_1^*)-(k_2-k_1)\frac{\partial}{\partial \theta} L(T_{k_2+1,n},\theta_1^*)
\Big]+o_P(1).
\end{multline}
By going along similar lines, we have 
\begin{multline}\label{part2_C_n,k1,k2}
\frac{k_1(k_2-k_1)}{n^{3/2}} J_{\theta^*_1}\left(\widehat{\theta}(T_{1,k_1})-\widehat{\theta}(T_{k_1+1,k_2})\right)=\\
\frac{1}{n^{3/2}}
\Big[
(k_2-k_1)\frac{\partial}{\partial \theta} L(T_{1,k_1},\theta_1^*)-k_1\frac{\partial}{\partial \theta} L(T_{k_1+1,k_2},\theta_1^*)
\Big]+o_P(1).
\end{multline}
Combining (\ref{part1_C_n,k1,k2}) and (\ref{part2_C_n,k1,k2}), we get
\begin{align*}
J_{\theta^*_1} C_{n,k_1,k_2}
&=
\frac{1}{n^{3/2}}
\Big[(n-(k_2-k_1))\frac{\partial}{\partial \theta} L(T_{k_1+1,k_2},\theta_1^*)
-(k_2-k_1) \big(\frac{\partial}{\partial \theta} L(T_{k_2+1,n},\theta_1^*)+ \frac{\partial}{\partial \theta} L(T_{1,k_1},\theta_1^*)\big) 
\Big] +o_P(1).  \\
&=
\frac{1}{\sqrt{n}}
\Big[
\frac{\partial}{\partial \theta} L(T_{k_1+1,k_2},\theta_1^*)
-\frac{(k_2-k_1)}{n}L(T_{1,n}) 
\Big] +o_P(1)\\
&=
\frac{1}{\sqrt{n}}
\Big[
\frac{\partial}{\partial \theta} L(T_{1,k_2},\theta_1^*)- \frac{\partial}{\partial \theta} L(T_{1,k_1},\theta_1^*)
-\frac{(k_2-k_1)}{n}L(T_{1,n}) 
\Big] +o_P(1) \\
&=
\frac{1}{\sqrt{n}}
\Big[
\big(\frac{\partial}{\partial \theta} L(T_{1,k_2},\theta_1^*) -\frac{k_2}{n}L(T_{1,n})\big)- 
\big(\frac{\partial}{\partial \theta} L(T_{1,k_1},\theta_1^*)
-\frac{k_1}{n}L(T_{1,n})\big) 
\Big] +o_P(1);
\end{align*}
i.e.,
\begin{equation}\label{Eq_b}
I^{-1/2}_{\theta^*_1}J_{\theta^*_1} C_{n,k_1,k_2}
=
\frac{I^{-1/2}_{\theta^*_1}}{\sqrt{n}}
\Big[
\big(\frac{\partial}{\partial \theta} L(T_{1,k_2},\theta_1^*) -\frac{k_2}{n}L(T_{1,n})\big)- 
\big(\frac{\partial}{\partial \theta} L(T_{1,k_1},\theta_1^*)
-\frac{k_1}{n}L(T_{1,n})\big) 
\Big] +o_P(1).
\end{equation}
Recall that, for any $\tau \in ]0,1]$,   
\[
\frac{\partial}{\partial \theta}L(T_{1,[n \tau]},\theta_1^*)=\sum_{t=1}^{[n \tau]}\frac{\partial}{\partial \theta}\ell_t(\theta_1^*),~
\text{ where } [n \tau] \text{ is the integer part of }  n \tau .
\] 
%
From  Lemma \ref{Lem00}(b) (applied to $\theta^*_1$), applying the central limit theorem for the martingale difference sequence $\left(\frac{\partial}{\partial \theta} \ell_t(\theta_1^*),\mathcal{F}_{t}\right)_{t \in \mathbb{Z}}$ (see \nocite{Billingsley1968} Billingsley (1968)), we have
\begin{align*}
\frac{1}{\sqrt{n}}\Big(\frac{\partial}{\partial \theta}L(T_{1,[n \tau_1]},\theta_1^*)-\frac{[n \tau_1]}{n}\frac{\partial}{\partial \theta}L(T_{1,n },\theta_1^*)\Big)
&=
\frac{1}{\sqrt{n}}\Big(\sum_{t=1}^{[n \tau_1]}\frac{\partial}{\partial \theta}\ell_t(\theta_1^*)-\frac{[n \tau_1]}{n}\sum_{t=1}^{n }\frac{\partial}{\partial \theta}\ell_t(\theta_1^*)\Big)\\
&~~
\limiteloin B_{I_{\theta^*_1}}(\tau_1)-\tau_1 B_{I_{\theta^*_1}}(1),
\end{align*}
where $B_{I_{\theta^*_1}}$ is a Gaussian process with covariance matrix $\min(s,t)I_{\theta^*_1}$. 
Hence, 
\begin{align*}
\frac{I^{-1/2}_{\theta^*_1}}{\sqrt{n}}\Big(\frac{\partial}{\partial \theta}L(T_{1,[n \tau_1]},\theta_1^*)-\frac{[n \tau_1]}{n}\frac{\partial}{\partial \theta}L(T_{1,n },\theta_1^*)\Big)
\limiteloin 
B_{d}(\tau_1)-\tau_1 B_{d}(1)=W_d(\tau_1) 
\end{align*}
 in  $D([0,1])$, where $B_d$ is a $d$-dimensional standard motion, and $W_d$ is a $d$-dimensional Brownian bridge.\\
 Similarly, we get 
  \begin{align*}
\frac{I^{-1/2}_{\theta^*_1}}{\sqrt{n}}\Big(\frac{\partial}{\partial \theta}L(T_{1,[n \tau_2]},\theta_1^*)-\frac{[n \tau_2]}{n}\frac{\partial}{\partial \theta}L(T_{1,n },\theta_1^*)\Big)
\limiteloin 
B_{d}(\tau_2)-\tau_2 B_{d}(1)=W_d(\tau_2) .
\end{align*}
Thus, as $n\rightarrow \infty$, it comes from (\ref{Eq_b}) that 
\begin{align*}
Q_{n,[n\tau_1],[n\tau_2]}&=
\big\| I^{-1/2}_{\theta^*_1}J_{\theta^*_1} C_{n,[n\tau_1],[n\tau_2]}\big\|^2
\limiteloin 
\left\|W_d(\tau_1)-W_d(\tau_2)\right\|^2~\text{ in }~ D([0,1]).
\end{align*}
Hence, for $n$ large enough, we have
\[
Q_n=\max_{\underset{k_1 <k_2-v_n}{v_n \leq k_1< k_2 \leq n-v_n}} Q_{n,k_1,k_2}=
\sup_{\frac{v_n}{n}\leq \tau_1 < \tau_2 \leq 1-\frac{v_n}{n}} Q_{n,[n\tau_1],[n\tau_2]} \limiteloin 
\sup_{0\leq \tau_1 < \tau_2 \leq 1} \left\|W_d(\tau_1)-W_d(\tau_2)\right\|^2\]
in~ D([0,1]). We conclude the proof of the theorem from Lemma \ref{Lem1}.
\begin{flushright}
$\blacksquare$ 
\end{flushright}

\subsubsection{Proof of Theorem \ref{th2}}
Recall that, under the alternative, the trajectory $(Y_1,\cdots,Y_n)$ satisfies
\begin{equation} \label{Eq_H1}
Y_{t}=\left\{
\begin{array}{ll}
Y^{(1)}_{t}~~\textrm{for}~~t \in [1,t^*_1]\cup[t^*_2+1,n],\\
\\
Y^{(2)}_{t}~~\textrm{for}~~t \in [t^*_1+1,t^*_2],\\
\end{array}
\right.
\end{equation} \\
where $(t^*_1,t^*_2)=([\tau^*_1 n],[\tau^*_2 n])$ (with $0<\tau^*_1<\tau^*_2<1$) and $\{Y^{(j)}_{t}, t \in \mathbb{Z}\}$ ($j=1,2$) is a stationary and ergodic solution of the model (\ref{Model}) depending on $\theta^{*}_j$ with $\theta^{*}_1 \neq \theta^{*}_2$.\\
We have $ \widehat{Q}_n = \underset{(k_1,k_2) \in \mathcal T_n}{\max}\widehat{Q}_{n,k_1,k_2} \geq \widehat{Q}_{n,t^*_1,t^*_2} $.  Then, it suffices to show that $\widehat{Q}_{n,t^*_1,t^*_2} \limiteproban +\infty$ 
to establish the theorem. 
For any $n \in \N$,  
\begin{flalign*}
\widehat{Q}_{n,t^*_1,t^*_2}&=C_{n,t^*_1,t^*_2}' \widehat{\Sigma}(u_n) C_{n,t^*_1,t^*_2}\\
\text{with}\hspace{1.5cm}&\\
 C_{n,t^*_1,t^*_2}&= \frac{t^*_2-t^*_1}{n^{3/2}} \left[ \big(n-(t^*_2-t^*_1)\big)\widehat{\theta}(T_{t^*_1+1,t^*_2}) - \left( t^*_1 \widehat{\theta}(T_{1,t^*_1}) + (n-t^*_2)\widehat{\theta}(T_{t^*_2+1,n})\right) \right]
&&
\end{flalign*}
and
\begin{multline*}\label{Sigma_un}
\widehat{\Sigma}(u_n)=\frac{1}{3} 
\big[
\widehat J(T_{1,u_n})  \widehat I(T_{1,u_n})^{-1}   \widehat J(T_{1,u_n}) 
+
\widehat J(T_{u_n+1,n-u_n})  \widehat I(T_{u_n+1,n-u_n})^{-1}  \widehat J(T_{u_n+1,n-u_n}) \\
+
\widehat J(T_{n-u_n+1,n})  \widehat I(T_{n-u_n+1,n})^{-1}  \widehat J(T_{n-u_n+1,n})
\big].
\end{multline*}
Moreover, for $n$ large enough, $\widehat{\theta}(T_{1,t^*_1}) \equalpsn \widehat{\theta}(T_{t^*_2+1,n})$ (from the consistency of the Poisson QMLE).
 Consequently, $C_{n,t^*_1,t^*_2}$ becomes 
\begin{equation}\label{C_(n,t^*_1,t^*_2)}
C_{n,t^*_1,t^*_2}= \frac{(t^*_2-t^*_1)(n-(t^*_2-t^*_1))}{n^{3/2}} \left( \widehat{\theta}(T_{1,t^*_1}) -\widehat{\theta}(T_{t^*_1+1,t^*_2})\right).
\end{equation}
Furthermore, by definition, the three matrices in the formula of $\widehat{\Sigma}_n(u_n)$ are positive semi-definite, 
and the first and the last one converge a.s. to same matrix which is positive definite. \\
Then, for $n$ large enough, we can  write
\begin{align*}
\widehat{Q}_{n,t^*_1,t^*_2}& \geq \\
& 
\frac{(t^*_2-t^*_1)^2(n-(t^*_2-t^*_1))^2}{n^{3}} \left( \widehat{\theta}(T_{1,t^*_1}) -\widehat{\theta}(T_{t^*_1+1,t^*_2})\right)'\\
& \hspace{2cm} \times \Big[
\widehat J(T_{1,u_n})  \widehat I(T_{1,u_n})^{-1}   \widehat J(T_{1,u_n}) +
\widehat J(T_{n-u_n+1,n})  \widehat I(T_{n-u_n+1,n})^{-1}  \widehat J(T_{n-u_n+1,n})
\Big]\\
& \hspace{2.5cm} \times  \left( \widehat{\theta}(T_{1,t^*_1}) -\widehat{\theta}(T_{t^*_1+1,t^*_2})\right)\\
& \geq 
R_n(\tau^*_1,\tau^*_2) \left( \widehat{\theta}(T_{1,t^*_1}) -\widehat{\theta}(T_{t^*_1+1,t^*_2})\right)'\\
& \hspace{2cm} \times \Big[
\widehat J(T_{1,u_n})  \widehat I(T_{1,u_n})^{-1}   \widehat J(T_{1,u_n}) +
\widehat J(T_{n-u_n+1,n})  \widehat I(T_{n-u_n+1,n})^{-1}  \widehat J(T_{n-u_n+1,n})
\Big]\\
& \hspace{2.5cm} \times  \left( \widehat{\theta}(T_{1,t^*_1}) -\widehat{\theta}(T_{t^*_1+1,t^*_2})\right)
\end{align*}
with
\[
R_n(\tau^*_1,\tau^*_2) = \frac{1}{n^3} \big[ \big(n(\tau^*_2-\tau^*_1)-1 \big) \big(n(1-(\tau^*_2-\tau^*_1)) -1 \big)\big]^2.
\]
From the asymptotic properties of the Poisson QMLE, we have 
\begin{align*}
\bullet&~~ \widehat{\theta}(T_{1,t^{*}})-\widehat{\theta}(T_{t^*_1+1,t^*_2}) \limitepsn \theta^{*}_{1}-\theta^{*}_{2}\neq 0,~~~\widehat{\theta}(T_{1,u_n}) \limitepsn \theta^{*}_1,
 ~~~\widehat{\theta}(T_{n-u_n+1,n}) \limitepsn \theta^{*}_2;\\
\bullet&~~\widehat J(T_{1,u_n})  \widehat I(T_{1,u_n})^{-1}   \widehat J(T_{1,u_n}) \limitepsn \Sigma^{(1)},~~~~
\widehat J(T_{n-u_n+1,n})  \widehat I(T_{n-u_n+1,n})^{-1}  \widehat J(T_{n-u_n+1,n}) \limitepsn \Sigma^{(1)},
\end{align*}
where
\begin{equation*}
\Sigma^{(1)}= J_1  I^{-1}_1 J_1 ~\text{ with } ~
J_1 =\E \Big[ \frac{1}{\lambda_{0}(\theta^*_1 )}  \frac{\partial \lambda_{0}(\theta^*_1 )}{ \partial \theta} \frac{\partial \lambda_{0}(\theta^*_1 )}{ \partial \theta'}  \Big] 
~\text{ and } ~
        I_1 =\E \Big[ \Big(\frac{Y_{0}}{\lambda_{0}(\theta^*_1 )}-1\Big)^2 \frac{\partial \lambda_{0}(\theta^*_1 )}{ \partial \theta} \frac{\partial \lambda_{0}(\theta^*_1 )}{ \partial \theta'}  \Big]. 
\end{equation*}
Therefore, since $\Sigma^{(1)}$ is positive definite and
$R_n(\tau^*_1,\tau^*_2) \limiten + \infty$, 
we deduce that $\widehat{Q}_{n,t^*_1,t^*_2} \limitepsn +\infty$. 
This completes the proof of the theorem.
\begin{flushright}
$\blacksquare$
\end{flushright}




\begin{thebibliography}{10}

\bibitem{Ahmad2016}
{\sc Ahmad, A., and Francq, C.}
\newblock Poisson qmle of count time series models.
\newblock {\em Journal of Time Series Analysis 37}, 3 (2016), 291--314.

\bibitem{Aston2012detecting}
{\sc Aston, J.~A., and Kirch, C.}
\newblock Detecting and estimating changes in dependent functional data.
\newblock {\em Journal of Multivariate Analysis 109\/} (2012a), 204--220.

\bibitem{Aston2012Evaluating}
{\sc Aston, J.~A., and Kirch, C.}
\newblock Evaluating stationarity via change-point alternatives with
  applications to fmri data.
\newblock {\em The Annals of Applied Statistics\/} (2012b), 1906--1948.

\bibitem{Billingsley1968}
{\sc Billingsley, P.}
\newblock {\em Convergence of probability Measures}.
\newblock John Wiley \& Sons, 1968.

\bibitem{Bucchia2014}
{\sc Bucchia, B.}
\newblock Testing for epidemic changes in the mean of a multiparameter
  stochastic process.
\newblock {\em Journal of Statistical Planning and Inference 150\/} (2014),
  124--141.

\bibitem{Csorgo1997}
{\sc Cs{\"o}rg{\"o}, M., and Horv{\'a}th, L.}
\newblock {\em Limit theorems in change-point analysis}.
\newblock Wiley New York, 1997.

\bibitem{Diop2020a}
{\sc Diop, M.~L., and Kengne, W.}
\newblock Consistent model selection procedure for general integer-valued time
  series.
\newblock {\em arXiv preprint arXiv:2002.08789\/} (2020a).

\bibitem{Diop2020b}
{\sc Diop, M.~L., and Kengne, W.}
\newblock Poisson \textsc{QMLE} for change-point detection in general
  integer-valued time series models.
\newblock {\em arXiv preprint arXiv:2007.13858\/} (2020b).

\bibitem{Diop2021}
{\sc Diop, M.~L., and Kengne, W.}
\newblock Piecewise autoregression for general integer-valued time series.
\newblock {\em Journal of Statistical Planning and Inference 211\/} (2021),
  271--286.

\bibitem{Doukhan2015}
{\sc Doukhan, P., and Kengne, W.}
\newblock Inference and testing for structural change in general poisson
  autoregressive models.
\newblock {\em Electronic Journal of Statistics 9\/} (2015), 1267--1314.

\bibitem{Ferland2006}
{\sc Ferland, R., Latour, A., and Oraichi, D.}
\newblock Integer-valued garch process.
\newblock {\em Journal of Time Series Analysis 27}, 6 (2006), 923--942.

\bibitem{Fokianos2009}
{\sc Fokianos, K., Rahbek, A., and Tj{\o}stheim, D.}
\newblock Poisson autoregression.
\newblock {\em Journal of the American Statistical Association 104}, 488
  (2009), 1430--1439.

\bibitem{Graiche2016}
{\sc Graiche, F., Merabet, D., and Hamadouche, D.}
\newblock Testing change in the variance with epidemic alternatives.
\newblock {\em Communications in Statistics-Theory and Methods 45}, 13 (2016),
  3822--3837.

\bibitem{Guan2007}
{\sc Guan, Z.}
\newblock Semiparametric tests for change-points with epidemic alternatives.
\newblock {\em Journal of statistical planning and inference 137}, 6 (2007),
  1748--1764.

\bibitem{Jaruvskova2011}
{\sc Jaru{\v{s}}kov{\'a}, D., and Piterbarg, V.~I.}
\newblock Log-likelihood ratio test for detecting transient change.
\newblock {\em Statistics \& probability letters 81}, 5 (2011), 552--559.

\bibitem{Levin1985}
{\sc Levin, B., and Kline, J.}
\newblock The cusum test of homogeneity with an application in spontaneous
  abortion epidemiology.
\newblock {\em Statistics in Medicine 4}, 4 (1985), 469--488.

\bibitem{Ravckauskas2004}
{\sc Ra{\v{c}}kauskas, A., and Suquet, C.}
\newblock H{\"o}lder norm test statistics for epidemic change.
\newblock {\em Journal of statistical planning and inference 126}, 2 (2004),
  495--520.

\bibitem{Ravckauskas2006}
{\sc Ra{\v{c}}kauskas, A., and Suquet, C.}
\newblock Testing epidemic changes of infinite dimensional parameters.
\newblock {\em Statistical Inference for Stochastic Processes 9}, 2 (2006),
  111--134.

\bibitem{Ramanayake2003}
{\sc Ramanayake, A., and Gupta, A.~K.}
\newblock Tests for an epidemic change in a sequence of exponentially
  distributed random variables.
\newblock {\em Biometrical Journal: Journal of Mathematical Methods in
  Biosciences 45}, 8 (2003), 946--958.

\bibitem{Souza2018}
{\sc Souza, J.~B., Reisen, V.~A., Franco, G.~C., Isp{\'a}ny, M., Bondon, P.,
  and Santos, J.~M.}
\newblock Generalized additive models with principal component analysis: an
  application to time series of respiratory disease and air pollution data.
\newblock {\em Journal of the Royal Statistical Society: Series C (Applied
  Statistics) 67}, 2 (2018), 453--480.

\bibitem{Weiss2008}
{\sc Wei{\ss}, C.~H.}
\newblock Thinning operations for modeling time series of counts-a survey.
\newblock {\em AStA Advances in Statistical Analysis 92}, 3 (2008), 319--341.

\bibitem{Weiss2019}
{\sc Wei{\ss}, C.~H., Feld, M. H.-J., Mamode~Khan, N., and Sunecher, Y.}
\newblock Inarma modeling of count time series.
\newblock {\em Stats 2}, 2 (2019), 284--320.

\bibitem{Weiss2014}
{\sc Wei{\ss}, C.~H., and Pollett, P.~K.}
\newblock Binomial autoregressive processes with density-dependent thinning.
\newblock {\em Journal of Time Series Analysis 35}, 2 (2014), 115--132.

\bibitem{Yao1993}
{\sc Yao, Q.}
\newblock Tests for change-points with epidemic alternatives.
\newblock {\em Biometrika 80}, 1 (1993), 179--191.

\bibitem{Zhu2011}
{\sc Zhu, F.}
\newblock A negative binomial integer-valued garch model.
\newblock {\em Journal of Time Series Analysis 32}, 1 (2011), 54--67.

\end{thebibliography}

 \end{document}